\documentclass[a4paper, 10pt]{amsart}
\usepackage{amsthm}
\usepackage{amssymb}
\usepackage{enumerate}
\usepackage{tabularx}
\usepackage[]{color}
\usepackage[left=2.2cm, right=2.2cm, top=2.5cm, bottom=2.5cm]{geometry}
\usepackage[colorlinks]{hyperref}
\usepackage{tikz}
\usepackage{multirow}
\usepackage{subcaption}
\usepackage{algorithm}
\usepackage{algorithmic}
\usepackage{xcolor}

\newcommand{\bm}[1]{\boldsymbol{#1}}
\newcommand{\bmr}[1]{\bm{\mr{#1}}}
\newcommand{\lj}{[ \hspace{-2pt} [}
\newcommand{\rj}{] \hspace{-2pt} ]}
\newcommand{\mb}[1]{\mathbb{#1}}
\newcommand{\mc}[1]{\mathcal{#1}}
\newcommand{\mr}[1]{\mathrm{#1}}
\newcommand{\jump}[1]{\lj #1 \rj}
\newcommand{\aver}[1]{ \{#1\}  }
\newcommand{\DGnorm}[1]{ \| #1\|_{\mr{DG}}}

\renewcommand{\d}[1]{\mathrm d \boldsymbol{#1}}
\newcommand\enorm[1]{|\!|\!| #1 |\!|\!|}
\def\curl{\ifmmode \mathrm{curl} \else \text{curl}\fi}
\def\div{\ifmmode \mathrm{div} \else \text{div}\fi}
\def\dim{\ifmmode \mathrm{dim} \else \text{dim}\fi}
\def\MTh{\mc{T}_h}
\def\MEh{\mc{E}_h}
\def\un{\bm{\mr n}}
\def\Ned{\ifmmode \text{N\'ed\'elec} \else \text{N\'ed\'elec} \fi}

\newcommand\comment[1]{}

\newcommand\substitute[2]{#2}

\newtheorem{assumption}{Assumption}
\newtheorem{theorem}{Theorem}
\newtheorem{lemma}{Lemma}

\definecolor{orange}{rgb}{1, 0.5, 0}

\allowdisplaybreaks[4]

\title[Quad-Curl Problem]{A Finite Element Method by Patch Reconstruction
for the Quad-Curl Problem Using Mixed Formulations}

\author[R. Li]{Ruo Li} \address{CAPT, LMAM and School of Mathematical
Sciences, Peking University, Beijing 100871, P.R. China}
\email{rli@math.pku.edu.cn}

\author[Q.-C. Liu]{Qicheng Liu} \address{School of Mathematical
Sciences, Peking University, Beijing 100871, P.R. China}
\email{qcliu@pku.edu.cn}

\author[S.-H. Zhao]{Shuhai Zhao} \address{School of Mathematical
Sciences, Peking University, Beijing 100871, P.R. China}
\email{shuhai@pku.edu.cn}

\begin{document}
\maketitle

\begin{abstract}
We develop a high order reconstructed discontinuous approximation
(RDA) method for solving a mixed formulation of the quad-curl problem 
in two and three dimensions. This mixed formulation is established 
by adding an auxiliary variable to control the divergence of the 
field. The approximation space for the original variables is
constructed by patch reconstruction with exactly one degree of freedom
per element in each dimension and the auxiliary variable is
approximated by the piecewise constant space. We prove the optimal 
convergence rate under the energy norm and also suboptimal $L^2$ 
convergence using a duality approach. Numerical results are provided 
to verify the theoretical analysis. \\ 
 \textbf{keywords}: quad-curl problem, mixed formulation, 
patch reconstruction
\end{abstract}

\section{Introduction}
\label{sec_introduction}
The quad-curl problem arises in many multiphysics
simulations, especially in inverse electromagnetic
scattering for inhomogeneous media ,magnetohydrodynamics and also 
Maxwell transmission eigenvalue problems. 
Therefore, it is important to design highly efficient and accurate 
numerical methods for quad-curl problems.

Finite element methods (FEMs) are a widely used numerical scheme for
solving partial differential equations. The presence of the quad-curl 
operator makes it difficult and challenging to design the conforming 
finite element space for quad-curl problems. We refer to \cite{Zhang2019curlcurl,
Hu2020simple,zhangqian2020} for somes recent works in constructing 
$H(\curl^2)$-conforming finite element spaces in two and three dimensions.
Due to the difficulties in discretizations of the quad-curl operator,
much attentions have been paid on using nonconforming elements 
such as \Ned's elements and completely discontinuous piecewise 
polynomials. We refer to \cite{zheng2011nonconforming, 
Hong2012discontinuous, Han2023hp} 
and the reference therein for works of this type. Another approach is 
devoted to mixed 
formulations. Some discussions can be found in \cite{Zhang2018mixed,Wang2019new,
Sun2016mixed}. Specifically, \cite{Zhang2018mixed} reduces the original 
problem to systems of low order equations by introducing intermediate variables, which makes 
the solution easier to approximate.


In this paper, we propose a mixed discontinuous Galerkin finite 
element method for the quad-curl problem with divergence-free variable.
A significant drawback of DG space is the large number of degrees of 
freedom in DG space, which results in high computational costs. This
drawback is a matter of concern. We follow the
methodology in \cite{Li2012efficient, Li2016discontinuous, 
Li2019reconstructed,Li2023curl} to apply
the patch reconstruction finite element method to the quad-curl problem. The
construction of the approximation space includes
creating an element patch for each element and solving a local least 
squares problem to obtain a polynomial basis function locally. 
Methods based on the reconstructed spaces are called reconstructed 
discontinuous approximation methods, which can approximate functions 
to high-order accuracy meanwhile 
inherits the flexibility on the mesh partition. One advantage of this
space is that it has very few degrees of freedom, which gives high 
approximation efficiency of finite element. 
The reconstructed space is a subspace of the standard DG space, 
so that we can borrow ideas from 
the interior penalty formulations to solve the quad-curl problem.
For the auxiliary variable, we use the piecewise constant space as the 
approximation space. Therefore, the mixed systems not grow much in 
size compared to the original system.
By adding penalty terms for 
both spaces, we do not need the two space to satisfy the discrete 
inf-sup condition.
We prove the convergence rates under
the energy norm and the $L^2$ norm, and numerical experiments
are conducted to verify the theoretical analysis and show that our 
algorithm is simple to implement and can reach high-order accuracy. 

The rest of this paper is organized as follows. In Section
\ref{sec_preliminaries}, we introduce the quad-curl problem with div-free condition
and give the basic notations about the Sobolev spaces and the
partition. In Section \ref{sec_space}, we introdude the RDA based finite element method. In
Section \ref{sec_problem}, we describe the mixed finite element method
for the quad-curl problem,
and prove that the convergence rate is optimal with respect to the energy norm
and suboptimal with respect to the $L^2$ norm. In Section
\ref{sec_numericalresults}, we carry out some numerical examples to
verify our theoretical results. A brief conclusion is given in Section
\ref{sec_conclusion}.


\section{preliminaries}
\label{sec_preliminaries}
Let $\Omega \subset \mb{R}^d (d = 2, 3)$ be a bounded polygonal
(polyhedral) domain with a Lipschitz boundary $\partial \Omega$. 

Given $\bm{f}\in H(\div^0,\Omega)$, we consider the following quad-curl problem
\begin{equation}
 \left\{
  \begin{aligned}
    \curl^{4}\bm{u}= \bm{f}, &\quad \text{ in } \Omega,\\
    \nabla \cdot\bm{u}=0, &\quad \text{ in } \Omega,
    \\
     \bm{u}\times\bm{n} = 0,  &\quad \text{ on } \partial \Omega,
    \\
     (\nabla \times \bm{u})\times\bm{n} = 0, &\quad \text{ on } \partial \Omega,
  \end{aligned}
 \right.
 \label{eq_quadcurl}
\end{equation}
We introduce an auxiliary variable $p$ to rewrite the problem as
\begin{equation}
  \left\{
   \begin{aligned}
     \curl^{4}\bm{u}+\nabla p= \bm{f}, &\quad \text{ in } \Omega,
     \\
     \nabla \cdot\bm{u}=0, &\quad \text{ in } \Omega,
     \\
      \bm{u}\times\bm{n} = 0,  &\quad \text{ on } \partial \Omega,
     \\
      (\nabla \times \bm{u})\times\bm{n} = 0, &\quad \text{ on } \partial \Omega,
     \\
      p = 0,  &\quad \text{on } \partial \Omega.
   \end{aligned}
  \right.
  \label{eq_quadcurl2}
 \end{equation}
 The weak form to the problem \eqref{eq_quadcurl2} is to find 
 $(\bm{u},p)\in H_0(\curl^2,\Omega)\times H_0^1(\Omega)$ such that 
  \begin{displaymath}
   \begin{aligned}
    (\curl^2 \bm{u}, \curl^2 \bm{v})+(\bm{v},\nabla p)&=(\bm{f},\bm{v}),\ \forall \bm{v}\in
     H_0(\curl^2,\Omega),\\
     -(\bm{u},\nabla q)&=0,\ \forall q\in H_0^1(\Omega),
   \end{aligned}
 \end{displaymath}
For the problem domain $\Omega$, we define
\begin{displaymath}
  \begin{aligned}
    H(\curl^s, \Omega) &:= 
      \left\{ \bm{v} \in L^2(\Omega)^d \ | \ \curl^j \bm{v} \in
      L^2(\Omega), 1\le j\le s \right\},  
  \end{aligned}
\end{displaymath}
\begin{displaymath}
  \begin{aligned}
    H_0(\curl^s, \Omega) &:= 
      \left\{ \bm{v} \in L^2(\Omega)^d \ | \ \curl^j \bm{v} \in
      L^2(\Omega), \ \curl^{j-1} \bm{v} = 0 \ \text{on}\  \partial \Omega,1\le j\le s \right\},  
  \end{aligned}
\end{displaymath}
 \begin{theorem}
 \cite{2019cc} Given $\bm{f}\in H(\div^0,\Omega)$, the problem \eqref{eq_quadcurl2} 
 admits a unique solution $(\bm{u},p)\in H_0(\curl^2,\Omega)\times H_0^1(\Omega)$ 
 with 
   \begin{displaymath}
    \|\bm{u}\|_{H(\curl^2,\Omega)} \le C\|\bm{f}\|,
 \end{displaymath}
 and $p=0$.
 \end{theorem}
 Let
\begin{displaymath}
  \begin{aligned}
    Y &:= \left\{ \bm{v} \in H_0(\curl, \Omega) \ | \ (\bm{v}, \nabla p)=0 , \forall p \in H_0^1(\Omega), \right\},  
  \end{aligned}
\end{displaymath}
\begin{lemma}
  \text{(Corollary 3.51 of \cite{Monk2003finite})} Suppose that $\Omega$ is a bounded Lipschitz domain. If $\Omega$ is simply connected and has a connected boundary, there is
  a $C>0$ such that for every $\bm{u}\in Y$
  \begin{equation}
      \|\bm{v}\|\le \|\curl\bm{v}\|
      \label{eq_zhimin3}
  \end{equation}
  \label{le_zhimin3}
\end{lemma}

\par Next, we define some notations about the mesh. Let $\MTh$ be a regular 
and quasi-uniform partition $\Omega$ into disjoint open triangles 
(tetrahedra). Let $\MEh$ denote the set of all $d-1$
dimensional faces of $\MTh$, and we decompose $\MEh$ into $\MEh = 
\MEh^i \cup \MEh^b$, where $\MEh^i$ and $\MEh^b$ are the 
sets of interior faces and boundary faces, respectively. We let 
\begin{displaymath}
  h_K := \text{diam}(K), \quad \forall K \in \MTh, \quad h_e :=
  \text{diam}(e), \quad \forall e \in \MEh, 
\end{displaymath}
and define $h := \max_{K \in \MTh} h_K$.
The quasi-uniformity of the mesh $\MTh$ is in the sense that there
exists a constant $\nu > 0$ such that $h \leq \nu \min_{K \in
\MTh} \rho_K$, where $\rho_K$ is the diameter of the
largest ball inscribed in $K$.


\section{Reconstructed Discontinuous Space}
\label{sec_space}
Now we introduce the local reconstruction operator to obtain the 
Reconstructed Discontinuous Approximation space. The first step is to
construct an element patch for each $K \in \MTh$.

For any element $K \in \MTh$ we construct an
element patch $S(K)$ which is an agglomeration of elements that
contain $K$ itself and some elements around $K$. There are a variety
of approaches to build the element patch and in this paper we
agglomerate elements to form the element patch recursively. For
element $K$, we first let $S_0(K) = \left\{ K \right\}$ and we define
$S_t(K)$ as
\begin{displaymath}
  S_t(K) = S_{t-1}(K) \cup \left\{ K'\ |\ \exists \widetilde{K} \in
  S_{t-1}(K)\ \text{s.t.}\ K' \cap \widetilde{K} = e \in \MEh\right\},
  \quad t = 1, 2, \cdots
\end{displaymath}
In the implementation of our code, at the depth $t$ we enlarge
$S_t(K)$ element by element and once $S_t(K)$ has collected
sufficiently large number of elements we stop the recursive procedure
and let $S(K) = S_t(K)$, otherwise we let $t = t + 1$ and continue the
recursion. The cardinality of $S(K)$ is denoted by $\# S(K)$.

We denote $\bm{x}_K$ the barycenter of
the element $K$ and mark barycenters of all elements as
\begin{displaymath}
  I(K) := \{ \bm{x}_{K'} \ |\ K' \in S(K) \},
\end{displaymath}

Let $U_h^0$ be the piecewise constant space, i.e.,
\begin{displaymath}
  U_h^0 := \{ v_h \in L^2(\Omega) \ | \  v_h|_K \in \mb{P}_0(K), \  
  \forall K \in \MTh\}.
\end{displaymath}
and $\bm{U}_h^0$ the $d$-dimensional piecewise constant space.

For any function $g\in \bm{U}_h^0$, we reconstruct a polynomial $\mc{R}_K g$
of degree $m$ on $S(K)$ by solving the least squares problem
\begin{equation} \label{eq:lsproblem}
  \mathcal R_Kg= \mathop{\arg \min}_{p\in\mb{P}_m(S(K))} \ \sum_{\bm{x}\in
    \mc{I}_K}\|\bm{g}(\bm{x})-\bm{p}(\bm{x})\|^2 \quad \text{s.t. } 
  \bm{g}(\bm{x}_K) = \bm{p}(\bm{x}_K), 
\end{equation}

The uniqueness condition for Problem~\eqref{eq:lsproblem} relates to
the location of the collocation points and $\# S(K)$. Following
\cite{Li2012efficient, Li2016discontinuous}, we make the following
assumption:
\begin{assumption}\label{as:unique}
For all $K\in\MTh$ and $p\in \mb{P}_m(S(K))$,
\begin{equation}
    p|_{\mc{I}_K}=0\quad\text{implies\quad} p|_{S(K)}\equiv0,
\end{equation}
\end{assumption}
The above assumption guarantees the uniqueness of the solution of
Problem~\eqref{eq:lsproblem} if $\# S(K)$ is greater than
$\text{dim}\,\mb{P}_m$. Hereafter, we assume that this assumption is
always valid.

The linear operator $\mc{R}$ can also be extended to act on smooth
functions in the following way. For any $\bm{g} \in H^{m+1}(\Omega)^d$, we define a
$\bm{g}_h \in \bmr{U}_h^0$ as 
\begin{displaymath}
  \bm{g}_h(\bm{x}_K) := \bm{g}(\bm{x}_K), \quad \forall K \in \MTh,
\end{displaymath}
and define $\mc{R} \bm{g} := \mc{R} \bm{g}_h$. Now we obtain the global
reconstruction operator $\mc{R}$.

Next, we will focus on the approximation properties of the operator
$\mc{R}$. We first define a constant for each element patch,
\substitute{}{
\begin{displaymath} 
    \Lambda(m,S(K)) := \max_{p \in
    \mb{P}_m(S(K))} \frac{\max_{\bm{x} \in S(K)} |p(\bm{x})|
    }{\max_{\bm{x} \in I(K)} |p(\bm{x})| },
\end{displaymath}}
and we refer to \cite{Li2016discontinuous, Li2023curl,Li2012efficient} for some 
discussions on the constants. Assumption \ref{as:unique} 
as well as the norm equivalence in finite dimensional spaces
 actually ensures $\Lambda(m,S(K))< \infty$.
\begin{lemma}
  For any element $K \in \MTh$, there holds
  \begin{equation}
    \substitute{}{\| \mc{R}_K \bm{g} \|_{L^\infty(K)} \leq (1+2
    \Lambda(m,S(K)) \sqrt{\# I(K)d}) \max_{\bm{x}\in I(K)} |\bm{g}|}
    , \quad \forall \bm{g} \in H^{m+1}(\Omega)^d. 
    \label{eq_stability}
  \end{equation}
  \label{le_stability}
\end{lemma}
\begin{proof}
Since $\bm{p} := \mc{R}\bm{g}$ is the solution to \eqref{eq:lsproblem}, for any
  $\varepsilon \in \mb{R}$ and 
  any $\bm{q} $ satisfying  the constraint in \eqref{eq:lsproblem}, 
  so does $\bm{p} + \varepsilon ( \bm{q} - \bm{g}(\bm{x}_{K})) $ and 
  \begin{displaymath}
    \sum_{\bm{x} \in I(K)} \| \bm{p}(\bm{x}) + \varepsilon (
    \bm{q}(\bm{x}) - \bm{g}(\bm{x}_{K}))- \bm{g}(\bm{x})
    \|_{l^2}^{2} \geq \sum_{\bm{x} \in I(K)} \| \bm{p}(\bm{x}) -
    \bm{g}(\bm{x}) \|_{l^2}^{2}.
  \end{displaymath}
  Since $\varepsilon$ is arbitrary, the above inequality implies 
  \begin{displaymath}
    \sum_{\bm{x} \in I(K)} (\bm{p}(\bm{x}) - \bm{g}(\bm{x})) \cdot
    (\bm{q} (\bm{x}) - \bm{g}(\bm{x}_{K})) = 0.
  \end{displaymath}
  By letting $\bm{q} = \bm{p}$, this orthogonal property indicates that

  \begin{equation}
    \sum_{\bm{x} \in I(K)} \| \bm{p}(\bm{x}) -
    \bm{g}(\bm{x}_{K}) \|_{l^2}^2 \leq  \sum_{\bm{x} \in I(K)} 
    \| \bm{g}(\bm{x}) - \bm{g}(\bm{x}_{K}) \|_{l^2}^2.
    \label{eq_pgh1}
  \end{equation}
  By definitions of the constants $ \Lambda(m,S(K))$, we get 
 
  \begin{displaymath}
    \begin{aligned}
      \| \bm{p}-\bm{g}(\bm{x}_K) \|_{L^\infty(K)}^2 &\leq \Lambda^2(m,S(K)) 
      \max_{\bm{x} \in I(K)} \| \bm{p}(\bm{x}) -
    \bm{g}(\bm{x}_{K}) \|_{l^2}^2 \leq \Lambda^2(m,S(K))\sum_{\bm{x} \in I(K)} 
    \| \bm{g}(\bm{x}) - \bm{g}(\bm{x}_{K}) \|_{l^2}^2\\
    &\le 4d \Lambda^2(m,S(K)) \# I(K) \max_{\bm{x}\in I(K)} |\bm{g}|^2,
    \end{aligned}
  \end{displaymath}
  hence
  \begin{displaymath}
    \begin{aligned}
      \| \bm{p} \|_{L^\infty(K)} &\leq \| \bm{p}-\bm{g}(\bm{x}_K) \|_{L^\infty(K)}
      +|\bm{g}(\bm{x}_K)|  &\le (1+2\Lambda(m,S(K)) \sqrt{\# I(K)d})\max_{\bm{x}\in I(K)} |\bm{g}|,
    \end{aligned}
  \end{displaymath}
  and completes the proof.
\end{proof}
\begin{assumption}
  For every element patch $S(K)(K \in \MTh)$,
  there exist constants $R$ and $r$ which are independent of $K$ such
  that $B_r \subset S(K) \subset B_R$, and $S(K)$ is star-shaped
  with respect to $B_r$, where $B_\rho$ is a disk with the radius
  $\rho$. 
  \label{as_patch}
\end{assumption}
From the stability result \eqref{eq_stability}, we can prove the
approximation results.  
\begin{lemma}
  For any $K \in \MTh $, there exists a constant $C$ such
  that 
  \begin{equation}
    \begin{aligned}
      \| \bm{g} - \mc{R} \bm{g} \|_{H^q(K)} &\leq C \Lambda_m h_K^{m
      + 1 - q} \| \bm{g} \|_{H^{m+1}(S(K))}, \quad 0 \leq q \leq m,
    \end{aligned}
    \label{eq_localapproximation}
  \end{equation}
  for any $\bm{g} \in H^{m+1}(\Omega)^d$, where we set 
  \begin{equation}
    \Lambda_m :=\max_{K \in \MTh}\left(1+\Lambda(m,S(K))\sqrt{\#I(K)d} \right).
    \label{eq_Lambdam}
  \end{equation}
  \label{le_localapproximation}
\end{lemma}
  We also present some useful lemmas commonly used in analyses concerning the curl operator.
  \begin{lemma}
    Let $\bm{v}_h\in \bmr{U}_h^m$, then there exists a $\bm{v}_h^c\in H_0(\curl,\Omega)$, such that
    \begin{equation}
      h^{-2}\|\bm{v}_h-\bm{v}_h^c\|^2+\|\curl_h(\bm{v}-\bm{v}_h^c)\|^2\lesssim \sum_{e \in \MEh}h_e^{-1}\|\jump{\bm{v}_h\times\un}\|_{L^2(e)}^2
      \label{eq_zhimin}
    \end{equation}
    \label{le_zhimin}
  \end{lemma}
  \begin{lemma}
    For any $\bm{v}_h\in \bmr{U}_h^m$, there exists a constant $C$ independent of the mesh size $h$, such that
    \begin{equation}
      \|\curl_h\bm{v}_h\|^2\lesssim \|\curl_h^2\bm{v}_h\|^2+\sum_{e \in \MEh}h_e^{-1}\|\jump{\curl\bm{v}_h\times\un}\|_{L^2(e)}^2+h_e^{-2}\|\jump{\bm{v}_h\times\un}\|_{L^2(e)}^2
      \label{eq_zhimin2}
    \end{equation}
    \label{le_zhimin2}
  \end{lemma}
  The above two lemmas can be found in \cite{zhimin2023}.

\section{Approximation to Quad-Curl Problem}
\label{sec_problem}
The mixed
discontinuous finite element method reads as follow: \textit{find
$\bm{u}_h \in \bm{\mr{U}}_h^m$ and $p \in U_h^0$, such that}
\begin{eqnarray}
  \begin{cases}
    &a(\bm{u}_h, \bm{v}_h) + b(p_h, \bm{v}_h) =
    F(\bm{v}_h) \quad \forall \bm{v}_h \in \bmr{U}_h^m, \\
    &b(q_h,\bm{u}_h) - c(p_h,q_h) = 0 \quad \forall q_h \in U_h^0. 
  \end{cases}
  \label{eq_mixform}
\end{eqnarray}
where
\begin{displaymath}
  \begin{aligned}
    a(\bm{u}_h, \bm{v}_h) &= \sum_{K \in \MTh} \int_{K} \curl ^{2}\bm{u} \cdot \curl ^{2}\bm{v}\
     \d{\bm{x}} + \sum_{e \in \MEh} \int_{e} 
    \left( \jump{\bm{u}\times\un} \cdot \aver{\curl^{3} \bm{v}} + 
    \jump{\curl \bm{u}\times\un} \cdot \aver{\curl^{2}\bm{v}} \right) \d{\bm{s}}\\
    &+ \sum_{e \in \MEh} \int_{e} 
    \left( \jump{\bm{v}\times\un} \cdot \aver{\curl^{3} \bm{u}} + 
    \jump{\curl \bm{v}\times\un} \cdot \aver{\curl^{2}\bm{u}} \right) \d{\bm{s}}
    \\
    &+ \sum_{e \in \MEh} \int_{e} \left(
    \mu_1 \jump{ \bm{u}\times\un} \cdot \jump{\bm{v}\times\un} + \mu_2 \jump{\curl \bm{u}\times\un}\cdot
    \jump{\curl \bm{v}\times\un} \right) \d{\bm{s}} \\
    &+ \sum_{K \in \MTh} \int_K (\nabla \cdot  \bm{u}_h)
    (\nabla \cdot  \bm{v}_h) \d{\bm{s}} + \sum_{e \in
    \MEh^I} \frac{\eta}{h_e} \int_e \jump{\un \cdot 
    \bm{u}_h} \jump{\un \cdot \bm{v}_h} \d{\bm{s}},
    \\
    b(p_h, \bm{v}_h) &= \sum_{K \in \MTh} \int_K (\nabla \cdot
     \bm{v}_h) p_h \d{\bm{x}} - \sum_{e \in \MEh^I}
    \int_e \jump{\un \cdot  \bm{v}_h} \aver{p_h}
    \d{\bm{s}}, \\
    c(p_h,q_h) &=  \sum_{e \in \MEh}  h_e \int_e \jump{p_h}
    \jump{q_h} \d{\bm{s}}, \\
    F(\bm{v}_h) &=  \sum_{K \in \MTh} \int_K \bm{f}\cdot\bm{v} \d{\bm{x}} .
  \end{aligned}
\end{displaymath}
The parameter $\mu_1$ and $\mu_2$ are positive penalties which are set
by
\begin{displaymath}
  \begin{aligned}
    \mu_1|_e = \frac{\eta}{h_e^3}, \quad \mu_2|_e =
    \frac{\eta}{h_e}. \quad \text{ for } e \in \MEh.
  \end{aligned}
\end{displaymath}
\substitute{}{
The global form of \eqref{eq_mixform} is defined by
\begin{equation}
  E(\bm{u}_h, p_h; \bm{v}_h,q_h) = F(\bm{v}_h), \quad
  \forall \bm{v}_h \in \bmr{U}_h^m, \quad q_h \in U_h^0,
\end{equation}
where
\begin{displaymath}
  E(\bm{u}, p; \bm{v},q) = a(\bm{u}, \bm{v})  + b(p, \bm{v}) - b(q, \bm{u}) + c(p,q),
\end{displaymath}
and the coefficient matrix of $E$ can be written as
\begin{displaymath}
  \begin{bmatrix}
    A & B^T\\ B & -C \\
  \end{bmatrix}
\end{displaymath}
where the matrices $A$, $B$ and $C$ associate with the bilinear
form $a(\cdot, \cdot)$, $b(\cdot, \cdot)$ and $c(\cdot, \cdot)$, 
respectively. We denote $n_e$ the number of elements in $\MTh$,
then
\begin{displaymath}
  A\in \mb{R}^{n_ed\times{n_ed}},\quad  B\in \mb{R}^{n_e\times n_ed},\quad  C\in \mb{R}^{n_e \times n_e}.
\end{displaymath}
}
We define the spaces $\bmr{V}_h := H(\curl^3, \Omega) \cap H(\div^0,
\Omega) + \bmr{U}_h^m$ and $Q_h := H_0^1(\Omega) + U_h^0$, \substitute{}{and
introduce the corresponding norms on $\bmr{U}_h^m$, $U_h^0$, and $\bmr{U}_h^m \times U_h^0$,
\begin{displaymath}
  \begin{aligned}
\DGnorm{\bm{v}}^2 &:= \sum_{K \in \MTh} \left(
   \| \curl^2 \bm{v} \|_{L^2(K)}^2 + 
  \| \nabla \cdot \bm{v} \|_{L^2(K)}^2 \right) \\
  &+ \sum_{e \in \MEh} h_e^{-3} \|
  \jump{ \bm{v}\times\un} \|_{L^2(e)}^2 + \sum_{e \in \MEh} h_e^{-1} \|
  \jump{ \curl\bm{v}\times\un} \|_{L^2(e)}^2+ \sum_{e \in \MEh^I} h_e^{-1} \|
  \jump{\un \cdot \bm{v}} \|_{L^2(e)}^2, \\
  \DGnorm{p}^2 &:=\sum_{e \in \MEh} h_e \|\jump{p} \|_{L^2(e)}^2,\\
  \DGnorm{(\bm{v}, p)} &:= \DGnorm{\bm{v}}+\DGnorm{p}.
\end{aligned}
\end{displaymath}}
Apparently $\DGnorm{\bm{v}}$ is a seminorm on $\bmr{U}_h^m$. We claim that the seminorm is a norm.
Given any $\DGnorm{\bm{v}}=0$, by definition 
    \begin{equation}
  \sum_{K \in \MTh}\| \curl^2 \bm{v} \|_{L^2(K)}=0;
  \label{1}
\end{equation}
  \begin{equation}
  \sum_{K \in \MTh}\| \nabla\cdot \bm{v} \|_{L^2(K)}=0;\,\sum_{e \in \MEh^I}\| \jump{ \bm{v}\cdot\un} \|_{L^2(e)}^2=0;
  \label{2}
\end{equation}
  \begin{equation}
  \sum_{e \in \MEh}  \| \jump{ \bm{v}\times\un} \|_{L^2(e)}^2=0;\,\, \sum_{e \in \MEh}  \| \jump{ \curl\bm{v}\times\un} \|_{L^2(e)}^2=0.
  \label{3}
\end{equation}

By Lemma \ref{le_zhimin} and \eqref{3}, $\bm{v}\in H_0(\curl,\Omega)$. Considering (\ref{2}), we can prove
\begin{displaymath}
  (\bm{v}, \nabla p)=-\sum_{K \in \MTh} \int_{K} p(\nabla\cdot\bm{v}) \d{\bm{x}}+\sum_{e \in \MEh^I}\int_e \aver{p} \jump{\bm{v}\cdot\un}=0 , \, \forall p \in H_0^1(\Omega),
\end{displaymath}
which means $\bm{v}\in Y$.
Lemma \ref{le_zhimin2} together with (\ref{1}) and (\ref{3}) gives $\| \curl \bm{v} \|_{L^2(\Omega)} =0$. 
Finally Lemma \ref{le_zhimin3} gives $\bm{v}=\bm{0}$, which means $\enorm{\bm{v}}$ is indeed a norm.

For the analyses we need another norm  
\substitute{}{
\begin{displaymath}
  \begin{aligned}
&\enorm{\bm{v}}^2 := \DGnorm{\bm{v}}^2 + \sum_{e \in
  \MEh}  \|h_e^{3/2} \aver{\curl^3 \bm{v}} \|_{L^2(e)}^2+ \sum_{e \in
  \MEh}  \|h_e^{1/2} \aver{\curl^2 \bm{v}} \|_{L^2(e)}^2,\\
  &\enorm{p} := \DGnorm{p},\quad  \enorm{(\bm{v}, p)} := \enorm{\bm{v}}+\enorm{p}.
  \end{aligned}
\end{displaymath}}
The norms $\DGnorm{\cdot}$ and $\enorm{\cdot}$ are equivalent
restricted on the space $\bmr{U}_h^m$. Obviously,
$\DGnorm{\cdot} \leq \enorm{\cdot} $. To verify $\enorm{\cdot} 
\leq C \DGnorm{\cdot}$,
for $e \in \MEh^i$, we denote $e=K^+\cap K^-$
  \begin{displaymath}
    \| h_e^{3/2} \aver{\curl^{3} \bm{u_h}} \|_{L^2(e)}
    \leq C (\| h_e^{3/2} \curl^{3} \bm{u_h}
    \|_{L^2(e \cap \partial K^+)} + \| h_e^{3/2} \curl^{3} \bm{u_h}
    \|_{L^2(e \cap \partial K^-)} ) .
  \end{displaymath}
  By the trace inequalities ,
  and the inverse inequality , we obtain that
  \begin{displaymath}
    \|h_e^{3/2} \curl^{3} \bm{u_h}\|_{L^2(e \cap \partial K^{\pm})}
    \leq 
      C \| \curl^{2} \bm{u_h} \|_{L^2(K^{\pm})}.
  \end{displaymath}
  For $e \in \MEh^b$, let $e \subset K$, and
  \begin{displaymath}
    \|h_e^{3/2} \curl^{3} \bm{u_h}\|_{L^2(e \cap \partial K)}
    \leq 
    C \| \curl^{2} \bm{u_h} \|_{L^2(K)}.
  \end{displaymath}
  The term $\| h_e^{1/2} \aver{\curl^{2} \bm{u_h}} \|_{L^2(e)}$ can 
  be bounded similarly. Thus, by summing over all $e \in \MEh$ ,
  we conclude that
  \begin{displaymath}
    \begin{aligned}
    \sum_{e \in \MEh} \| h_e^{1/2} \aver{\curl^{2} \bm{u_h}} \|_{L^2(e)}^2 
    + \sum_{e \in \MEh} \| h_e^{3/2} \aver{\curl^{3} \bm{u_h}}
    \|_{L^2(e)}^2  \leq C \sum_{K \in \MTh} \| \curl^{2} \bm{u_h} \|_{L^2(K)}^2.
    \end{aligned}
  \end{displaymath}

Next, we can prove the coercivity for the form $E$,
\begin{theorem}
  There exists a positive constant $C$ independent of $h$ such that,
  for all $(\bm{v}_h, q_h) \in \bmr{U}_h^m \times U_h^0$, 
  \begin{displaymath}
    E(\bm{v}_h, q_h; \bm{v}_h, q_h) \geq C \enorm{(\bm{v}_h,
    q_h)}^2 .
  \end{displaymath}
  \label{th_garding}
\end{theorem}
\begin{proof}
  We first note that
  \begin{displaymath}
    E(\bm{v}_h, q_h; \bm{v}_h, q_h) = a(\bm{v}_h, \bm{v}_h) + c(q_h,q_h).
  \end{displaymath}
  From the norm equivalence claimed above we only need to establish the coercivity of $a$
  over the norm $\DGnorm{\cdot}$. From Cauchy-Schwartz inequalities, trace inequalities
  , and inverse inequalities, 
  \begin{equation}
    -\sum_{e \in \MEh} \int_{e} 2 \jump{\bm{u_h}\times\un} \cdot
    \aver{\curl^{3}\bm{u_h}} \d{\bm{s}} \geq -\sum_{e \in
    \MEh^b} \frac{1}{\epsilon} \| h_e^{-3/2} \jump{\bm{u_h}\times\un}
    \|^2_{L^2(e)} - C\epsilon \sum_{K \in \MTh} \| \curl^{2} \bm{u_h}
    \|^2_{L^2(K)}.
    \label{eq_eb1}
  \end{equation}
also
  \begin{equation}
    -\sum_{e \in \MEh} \int_{e} 2 \jump{\curl \bm{u_h}\times\un}
    \aver{\curl^{2} \bm{u_h}} \d{\bm{s}} \geq -\frac{1}{\epsilon} \sum_{e
    \in \MEh} \| h_e^{-1/2} \jump{\curl \bm{u_h}\times\un} \|^2_{L^2(e)}
    - C\epsilon \sum_{K \in \MTh} \| \curl^{2}\bm{u_h} \|^2_{L^2(K)},
    \label{eq_eob2}
  \end{equation}
Therefore
  \begin{displaymath}
    \begin{aligned}
      a(\bm{u_h}, \bm{u_h}) &\geq (1 - C\epsilon) 
      \sum_{K \in \MTh} \|  \curl^{2} \bm{u_h}\|^2_{L^2(K)} +\sum_{K \in \MTh} \|  \nabla\cdot \bm{u_h}\|^2_{L^2(K)} \\
    &+ (\eta - \frac{1}{\epsilon}) \sum_{e
    \in \MEh} (\| h_e^{-1/2} \jump{ \curl \bm{u_h}\times\un} \|^2_{L^2(e)} 
    + \| h_e^{-3/2} \jump{\bm{u_h}\times\un} \|^2_{L^2(e)}) + \sum_{e
    \in \MEh^I}\| h_e^{-1/2} \jump{\bm{u_h}\cdot\un} \|^2_{L^2(e)}
    \end{aligned}
  \end{displaymath}
  for any $\epsilon > 0$. We can let $\epsilon = 1/(2C)$ and
  select a sufficiently large $\eta$ to ensure $a(\bm{u_h}, \bm{u_h}) \geq C
  \DGnorm{\bm{u_h}}^2$, which completes the proof.
\end{proof}

\begin{theorem}
  Let $(\bm{u}, p)$ be the solution of the quad-curl equations
  \eqref{eq_quadcurl2}, $(\bm{u}_h, p_h)$ be the solution of the
  approximation problem \eqref{eq_mixform}. Then for all $(\bm{v}_h,
  q_h) \in \bmr{U}_h^m \times U_h^0$, we have
  \begin{displaymath}
    E(\bm{u} - \bm{u}_h, p-p_h; \bm{v}_h, q_h) = 0.
  \end{displaymath}
  \label{th_orthogonality}
\end{theorem}
\begin{proof}
  \begin{displaymath}
    \begin{aligned}
      E(\bm{u} - \bm{u}_h, &p-p_h; \bm{v}_h, q_h) = \sum_{K \in
      \MTh} \int_K \curl^2 (\bm{u} - \bm{u}_h) \cdot\curl^2
      \bm{v}_h \d{\bm{x}}  \\
      &- \sum_{e \in \MEh} \int_e \jump{
      \bm{v}_h \times \un} \cdot \aver{ \curl^3  (\bm{u}-\bm{u}_h)}
      \d{\bm{s}} - \sum_{e \in \MEh} \int_e \jump{ \bm{u}_h 
      \times\un} \cdot \aver{\curl^3  \bm{v}_h} 
      \d{\bm{s}} \\
      &- \sum_{e \in \MEh} \int_e \jump{\curl
      \bm{v}_h \times \un} \cdot \aver{ \curl^2  (\bm{u}-\bm{u}_h)}
      \d{\bm{s}} - \sum_{e \in \MEh} \int_e \jump{ \curl\bm{u}_h 
      \times\un} \cdot \aver{\curl^2  \bm{v}_h} 
      \d{\bm{s}} \\
      &- \sum_{e \in \MEh}
      \mu_1 \int_e \jump{\bm{u}_h\times\un } \cdot 
      \jump{ \bm{v}_h\times\un } \d{\bm{s}}- \sum_{e \in \MEh}
      \mu_2 \int_e \jump{\curl\bm{u}_h\times\un } \cdot 
      \jump{\curl\bm{v}_h\times\un } \d{\bm{s}} \\
      &+\sum_{K \in \MTh} \int_K \nabla \cdot ( (\bm{u} -
      \bm{u}_h)) \nabla \cdot \bm{v}_h\d{\bm{x}}.
      -\sum_{e \in \MEh^I} \frac{\eta}{h_e} \int_e \jump{\un \cdot
       \bm{u}_h} \jump{\un \cdot 
       \bm{v}_h} \d{\bm{s}} \\
      &+b(p - p_h, \bm{v}_h) - b(q_h, \bm{u}-\bm{u}_h) + c(p-p_h,q_h)
      \\
      & \quad =\sum_{K \in \MTh} \int_K  \curl^2 \bm{u} \cdot
      \curl^2 \bm{v}_h \d{\bm{x}} + \sum_{K \in \MTh} \int_K
      (\nabla \cdot \bm{u})(\nabla \cdot \bm{v}_h )\d{\bm{x}} \\
      &- \sum_{e \in \MEh} \int_e
      \jump{ \bm{v}_h\times\un } \cdot \aver{ \curl^3
      \bm{u}} \d{\bm{s}} - \sum_{e \in \MEh} \int_e
      \jump{\curl \bm{v}_h\times\un} \cdot \aver{ \curl^2
      \bm{u}} \d{\bm{s}} \\
      &- \int_{\Omega} \bm{f} \cdot \bm{v}_h
      \d{\bm{x}} + b(p,
      \bm{v}_h) - b(q_h, \bm{u}) + c(p,q_h).
    \end{aligned}
  \end{displaymath}
  Since $p = 0$ and $\div \bm{u} = 0$, we have
  \begin{displaymath}
    b(p,\bm{v}_h) = b(q_h, \bm{u}) = c(p,q_h) = 0,
  \end{displaymath}
  which implies $E(\bm{u} - \bm{u}_h, p-p_h; \bm{v}_h, q_h) = 0$,
  and completes the proof.
\end{proof}
Before proving the error estimates, we need to establish the interpolation error estimate of the reconstruction
operator.
\begin{lemma}
  For $0 \leq h \leq h_0$ and $m \geq 2$, there exists a constant $C$ 
  such that
  \begin{equation}
    \substitute{}{\enorm{\bm{v} - \mc{R} \bm{v}} \leq C \Lambda_m h^{m-1} \| \bm{v}
    \|_{H^{m+1}(\Omega)}, \quad \forall \bm{v} \in
    H^{s}(\Omega), \quad s = \max(4, m+1).}
    \label{eq_interpolation_err}
  \end{equation}
  \label{le_interpolation_err}
\end{lemma}
\begin{proof}
  From Lemma
  \ref{le_localapproximation}, we can show that
  \begin{displaymath}
    \begin{aligned}
      \sum_{K \in \MTh} \| \curl^{2} \bm{v} - \curl^{2}(\mc{R} \bm{v}) \|^2_{L^2(K)} 
      &\leq \sum_{K \in \MTh} C \Lambda_m^2 h_K^{2m-2} \|  \bm{v}
      \|^2_{H^{m+1}(S(K))} \\
      &\leq C \Lambda_m^2 h^{2m-2} \|  \bm{v} \|^2_{H^{m+1}(\Omega)},
    \end{aligned}
  \end{displaymath}
  also
  \begin{displaymath}
    \begin{aligned}
      \sum_{K \in \MTh} \|  \nabla\cdot\bm{v} - \nabla\cdot(\mc{R} \bm{v}) \|^2_{L^2(K)} 
      &\leq C \Lambda_m^2 h^{2m} \| \bm{v} \|^2_{H^{m+1}(\Omega)}, \\
    \end{aligned}
  \end{displaymath}
  By the trace estimate and the mesh regularity,
  \begin{displaymath}
    \begin{aligned}
      \sum_{e \in \MEh}  h_e^{-1} \| \jump{\curl ( \bm{v} - \mc{R} \bm{v})\times\un} 
      \|^2_{L^2(e)} &\leq C \sum_{K \in \MTh} 
      \left( h_K^{-2} \| \curl( \bm{v} - \mc{R}\bm{v}) \|^2_{L^2(K)}
      + \| \curl ( \bm{v} - \mc{R}\bm{v}) \|^2_{H^1(K)} \right) \\
      &\leq C \Lambda_m^2 h^{2m-2} \|  \bm{v}
      \|^2_{H^{m+1}(\Omega)},
    \end{aligned}
  \end{displaymath}
  \begin{displaymath}
    \begin{aligned}
      \sum_{e \in \MEh}  h_e^{-3} \| \jump{  (\bm{v} - \mc{R} \bm{v})\times\un} 
      \|^2_{L^2(e)}&\leq C \Lambda_m^2 h^{2m-2} \|  \bm{v}
      \|^2_{H^{m+1}(\Omega)},
    \end{aligned}
  \end{displaymath}
  and also
  \begin{displaymath}
    \begin{aligned}
      \sum_{e \in \MEh^I}  h_e^{-1} \| \jump{  (\bm{v} - \mc{R} \bm{v})\cdot\un} 
      \|^2_{L^2(e)}
      &\leq C \Lambda_m^2 h^{2m} \|  \bm{v}
      \|^2_{H^{m+1}(\Omega)} .
    \end{aligned}
  \end{displaymath}
  The other terms can be estimated by trace estimates and interpolation error estimates similarly.
\end{proof}
\begin{theorem}
  \substitute{}{Let $(\bm{u}, p)$ be the solution of the quad-curl equations
  \eqref{eq_quadcurl2},}and suppose $\bm{u} \in H^{s}(\Omega)$, where $s = \max(4,m+1)$.
  For sufficient large $\eta > 0$, the error $(\bm{u} -
  \bm{u}_h, p-p_h)$ satisfies
  \begin{equation}
    \enorm{(\bm{u} - \bm{u}_h, p-p_h)} \leq  C\Lambda_m h^{m-1} \| \bm{u} \|_{H^{m+1}(\Omega)}.
    \label{eq_err1}
  \end{equation}
  \label{th_err1}
\end{theorem}
\begin{proof}
  \begin{displaymath}
    \begin{aligned}
      C\enorm{(\bm{u}_h-\mc{R}\bm{u},p_h)}^2&\le E(\bm{u}_h-\mc{R}\bm{u},p_h,\bm{u}_h-\mc{R}\bm{u},p_h)=E(\bm{u}-\mc{R}\bm{u},0;\bm{u}_h-\mc{R}\bm{u},p_h)\\
      &= a (\bm{u}-\mc{R}\bm{u},\bm{u}_h-\mc{R}\bm{u}) -b (p_h,\bm{u}-\mc{R}\bm{u}) \\
      &=:E_1+E_2.
    \end{aligned}
  \end{displaymath}
By Cauchy-Schwartz inequality, 
  \begin{displaymath}
    \begin{aligned}
  E_1   &\le \enorm{\bm{u} -\mc{R}\bm{u}}\enorm{\bm{u}_h - \mc{R}\bm{u}}\\
  &\le C\Lambda_m h^{m-1} \| \bm{u} \|_{H^{m+1}(\Omega)}\enorm{\bm{u}_h-\mc{R}\bm{u}}.\\\\
\end{aligned}
\end{displaymath}
  Then we estimate $E_2$,
\begin{displaymath}
  \begin{aligned}
  E_2 &= \sum_{K \in \MTh} \int_{\partial K} p_h(\bm{u}-\mc{R}\bm{u}) \cdot \bm{n} \d{\bm{s}}- \sum_{e \in \MEh^I}\int_e \aver{p_h} \jump{(\bm{u}-\mc{R}\bm{u})\cdot\un}\\
    &= \sum_{e \in \MEh}
    \int_e \jump{p_h} \aver{\bm{u}-\mc{R}\bm{u}}\\
    &\le \enorm{p_h}\left(\sum_{e \in \MEh} \int_e h^{-1}|\aver{\bm{u}-\mc{R}\bm{u}}|^2 \d{\bm{s}}\right)^{1/2}\\
    &\le C\Lambda_m h^{m} \| \bm{u} \|_{H^{m+1}(\Omega)}\enorm{p_h}.
      \end{aligned}
  \end{displaymath}
  Therefore
  \begin{displaymath}
    \enorm{(\bm{u}_h-\mc{R}\bm{u},p_h)}\leq  C\Lambda_m h^{m-1} \| \bm{u} \|_{H^{m+1}(\Omega)}.
  \end{displaymath}
  The proof is complished by 
  \begin{displaymath}
    \enorm{(\bm{u}_h-\bm{u}_h,p_h)}\leq  \enorm{\bm{u}-\mc{R}\bm{u}}+ \enorm{(\bm{u}_h-\mc{R}\bm{u},p_h)},
  \end{displaymath}
  using Lemma \ref{le_interpolation_err}.
\end{proof}
Now we turn to $L^2$ estimates. We introduce an auxiliary problem
\begin{equation}
\left\{
    \begin{aligned}
  \curl^4 \bm{w} + \nabla \xi &= \bm{u}-\bm{u}_h \quad \text{in } \Omega, \\
  \nabla \cdot \bm{w} &= 0, \quad \text{in } \Omega, \\
   \bm{w}\times\bm{n} &= 0, \quad \text{on } \partial \Omega, \\
    \nabla\times\bm{w} &= 0, \quad \text{on } \partial \Omega, \\
  \xi &= 0, \quad \text{on } \partial \Omega.
   \end{aligned}
  \right.
\label{eq_quadcurl3}
\end{equation}
We assume the regularity estimate as in \cite{wgcurl,Han2023hp}.
  \begin{equation}
    \| \bm{w} \|_{H^4(\Omega)}+\| \xi\|_{H^1(\Omega)} \le  C\|\bm{u} - \bm{u}_h\|.
    \label{eq_err1}
  \end{equation}
\begin{theorem}
  Under the same assumptions as Theorem \ref{th_err1}, the $L^2$ error satisfies
  \begin{equation}
    \|\bm{u} - \bm{u}_h\| \leq  C\Lambda_m h^{m} \| \bm{u} \|_{H^{m+1}(\Omega)}.
    \label{eq_err2}
  \end{equation}
  \label{th_err2}
\end{theorem}
\begin{proof}
  By taking inner product with respect to $\bm{u}-\bm{u}_h$ in the first equation of \eqref{eq_quadcurl3}
  \begin{displaymath}
    \begin{aligned}
      \|\bm{u}-\bm{u}_h\|^2 &= a (\bm{w},\bm{u} - \bm{u}_h) +b (\xi,\bm{u}-\bm{u}_h) \\
      &= a (\bm{w}-\mc{R}\bm{w},\bm{u} - \bm{u}_h) -b (p-p_h,\mc{R}\bm{w}) 
      +b(\xi-R\xi,\bm{u} - \bm{u}_h) +b(R\xi,\bm{u} - \bm{u}_h)\\
      &=:e_1+e_2+e_3+e_4,
    \end{aligned}
  \end{displaymath}
  where $R: H^1(\Omega) \rightarrow U_h^0$ is the local $L^2$ projection, satisfying
  \begin{displaymath}
    \begin{aligned}
  \|\xi-R\xi\|\le Ch\|\nabla \xi\|, \, \forall \xi \in H^1(\Omega).
    \end{aligned}
\end{displaymath}
  We estimate the four terms respectively. First
  \begin{displaymath}
    \begin{aligned}
  e_1 &\le C\enorm{\bm{w}-\mc{R}\bm{w}}\enorm{\bm{u} - \bm{u}_h}\\
  &\le C\Lambda_m h^{m+1} \| \bm{u} \|_{H^{m+1}(\Omega)}\| \bm{w} \|_{H^{4}(\Omega)}.
    \end{aligned}
\end{displaymath}
  Since $p = 0$ and $\div \bm{w} = 0$, we have
  \begin{displaymath}
      \begin{aligned}
    e_2 &= b(p_h,\mc{R}\bm{w}) = -b(p_h,\bm{w}-\mc{R}\bm{w}) \\
    &= \sum_{e \in \MEh}
    \int_e \jump{p_h} \aver{\bm{w}-\mc{R}\bm{w}}\d{\bm{s}}\\
    &\le \enorm{p_h}\left(\sum_{e \in \MEh} \int_e h^{-1}|\aver{\bm{w}-\mc{R}\bm{w}}|^2 \d{\bm{s}}\right)^{1/2}\\
    &\le C\Lambda_m h^{m+2} \| \bm{u} \|_{H^{m+1}(\Omega)}\| \bm{w} \|_{H^{4}(\Omega)}.
      \end{aligned}
  \end{displaymath}
  Also by definition,
  \begin{displaymath}
    \begin{aligned}
    e_3=& \sum_{K \in \MTh} \int_K (\nabla \cdot
    (\bm{u}-\bm{u}_h)) (\xi-R\xi) \d{\bm{x}} - \sum_{e \in \MEh^I}
   \int_e \jump{\un \cdot(\bm{u}-\bm{u}_h)} \aver{\xi-R\xi}
   \d{\bm{s}}\\
   &\le \|\xi-R\xi\|\|\div(\bm{u}-\bm{u}_h)\|+\left(\sum_{e \in \MEh^I} \int_e h|\aver{\xi-R\xi}|^2 \d{\bm{s}}\right)^{1/2}\enorm{\bm{u}-\bm{u}_h} \\
   & \le C\Lambda_m h^{m} \| \bm{u} \|_{H^{m+1}(\Omega)}\| \xi \|_{H^{1}(\Omega)}.
    \end{aligned}
  \end{displaymath}
  The error of the numerical solution satisfies
  \begin{displaymath}
    -b(q_h,\bm{u}-\bm{u}_h)+ \sum_{e \in \MEh} \int_e h\jump{p-p_h}\cdot\jump{q_h} \d{\bm{s}}=0, \, \forall q_h\in U_h^0.
  \end{displaymath}
  Using trace inequality, 
  \begin{displaymath}
    \begin{aligned}
    e_4&= -\sum_{e \in \MEh} \int_e h\jump{p-p_h}\cdot\jump{\xi-R\xi}\d{\bm{s}}
    \le\left(\sum_{e \in \MEh} \int_e h|\jump{p_h}|^2 \d{\bm{s}}\right)^{1/2}\left(\sum_{e \in \MEh} \int_e h|\jump{\xi-R\xi}|^2 \d{\bm{s}}\right)^{1/2}\\
     &\le h\|\xi\|_{H^1(\Omega)}\enorm{p_h}\le C\Lambda_m h^{m} \| \bm{u} \|_{H^{m+1}(\Omega)}\| \xi \|_{H^{1}(\Omega)}.
    \end{aligned}
  \end{displaymath}
  The proof is accomplished by combining the estimates of $e_1, e_2, e_3, e_4$ and the regularity assumption \eqref{eq_err1}.
\end{proof}

\section{Numerical Results}
\label{sec_numericalresults}
In this section, we perform numerical experiments to test the
performance of our method. We shall solve the following quad-curl problem with non-homogeneous
boundary conditions:
\begin{equation}
  \left\{
   \begin{aligned}
     \curl^{4}\bm{u}= \bm{f}, &\quad \text{ in } \Omega,\\
     \nabla \cdot\bm{u}=0, &\quad \text{ in } \Omega,
     \\
      \bm{u}\times\bm{n} = \bm{g}_1,  &\quad \text{ on } \partial \Omega,
     \\
      (\nabla \times \bm{u})\times\bm{n} = \bm{g}_2, &\quad \text{ on } \partial \Omega,
   \end{aligned}
  \right.
  \label{eq_quadcurl4}
 \end{equation}
 in such case the right hand side $F$ takes the form 
 \begin{equation}
  F(\bm{v}_h) =  \sum_{K \in \MTh} \int_K \bm{f}\cdot\bm{v} \d{\bm{x}} 
  + \sum_{e \in \MEh^b} \int_{e} \left( \bm{g}_1 \cdot(\curl^{3}\bm{v_h}
   + \mu_1\, \bm{v_h}\times \un) + \bm{g}_2 \cdot(\curl^{2}\bm{v_h} + \mu_2\, \curl\bm{v_h}\times\un) \right)
  \d{\bm{s}}.
 \end{equation}
In the test examples, the right hand side as well as boundary conditions are chosen according to the
exact solution. 

\begin{figure}[htbp]
  \centering
  \includegraphics[width=0.30\textwidth]{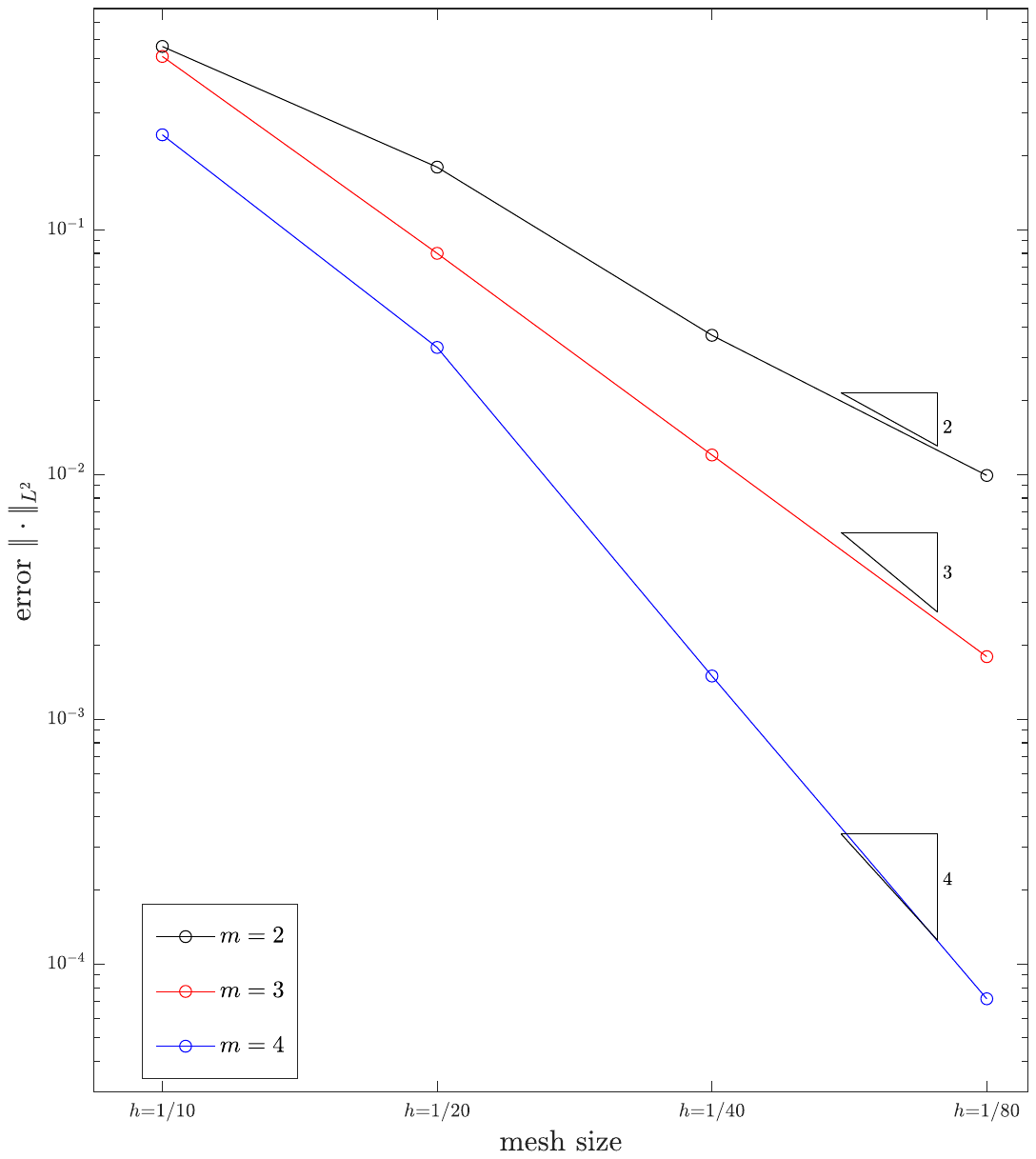}
  \hspace{30pt}
  \includegraphics[width=0.30\textwidth]{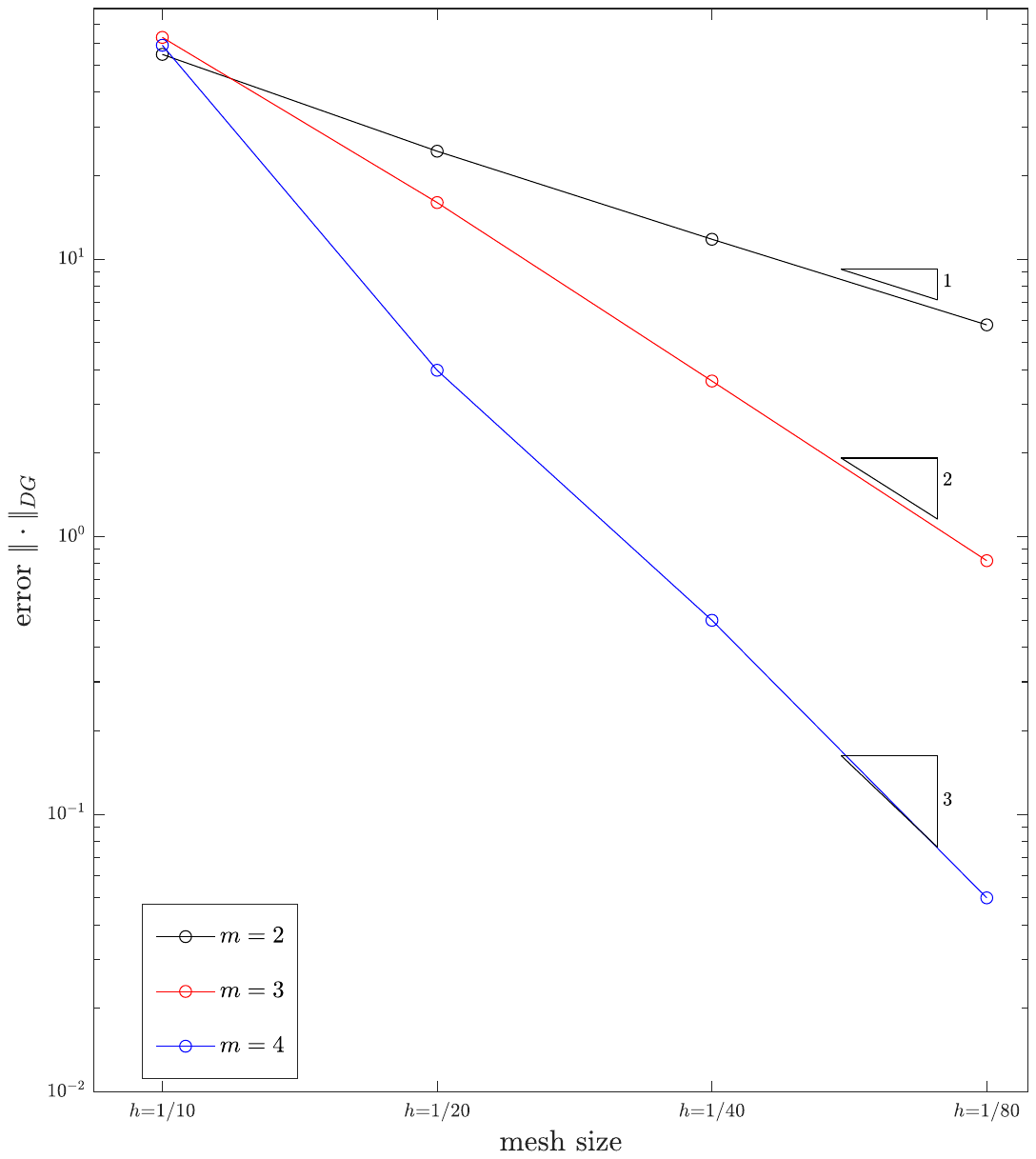}
  \caption{The convergence histories under the $\| \cdot \|_{L^2(\Omega)}$ 
  (left) and the$\DGnorm{\cdot}$ (right) in Example 1.}
  \label{fig_ex1err}
\end{figure}

\noindent \textbf{Example 1.} 
We first give an example on the 2D domain $\Omega=(0,1)^2$
\begin{displaymath}
  \bm{u}(x, y) = \begin{bmatrix}
    3\pi \sin^2(\pi y)\cos(\pi y)\sin^3(\pi x)\\
    -3\pi \sin^2(\pi x)\cos(\pi x)\sin^3(\pi y) \\
  \end{bmatrix},
\end{displaymath}
\begin{figure}
  \centering
  \includegraphics[width=0.3\textwidth]{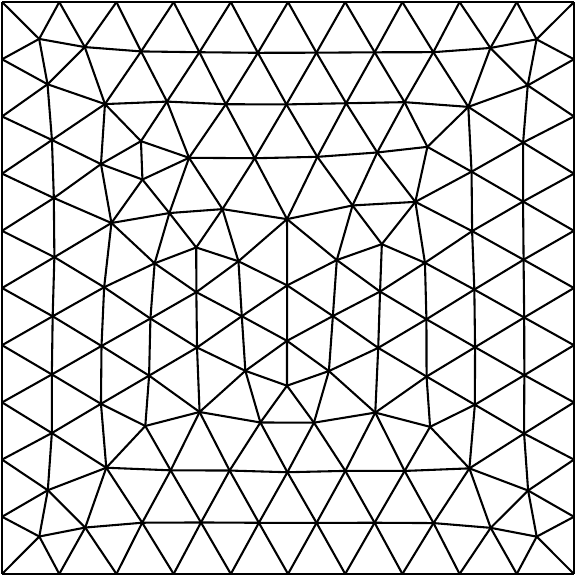}
  \hspace{25pt}
  \includegraphics[width=0.3\textwidth]{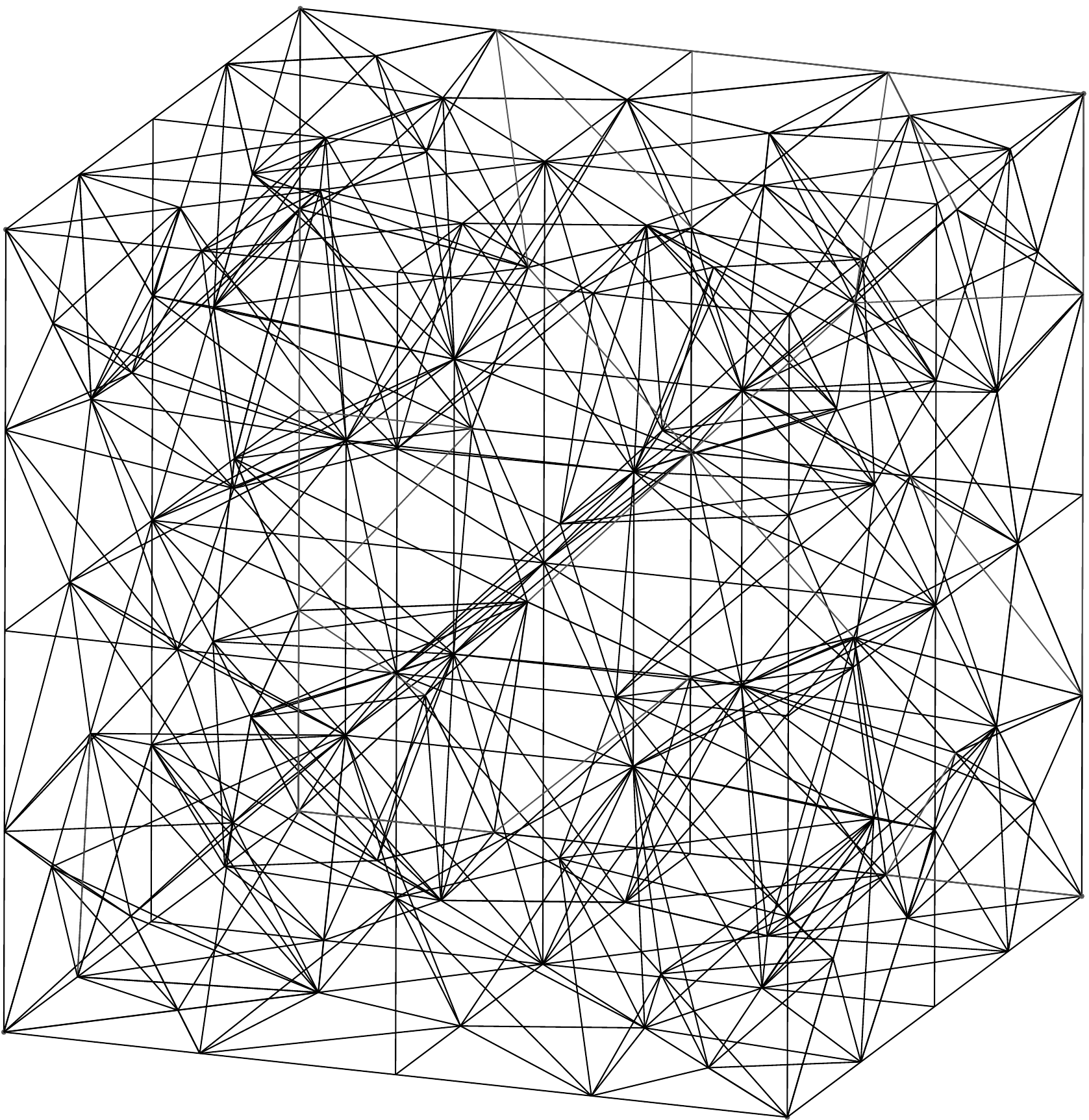}
  \caption{2d triangular partition with $ h = 1/10$ (left) / 3d
  tetrahedral partition with $h = 1/4$ (right).}
  \label{fig_partition}
\end{figure}
We solve the quad-curl problem on a sequence of meshes with
$h = 1/10, 1/20, 1/40, 1/80$. The convergence histories under 
the $\DGnorm{\cdot}$ and
$\| \cdot \|_{L^2(\Omega)}$ are shown in
Fig.~\ref{fig_ex1err}. We observe the optimal convergence of DG norm and suboptimal convergence of $L^2$ norm.

\begin{table}[htp]
\begin{minipage}[t]{0.3\textwidth}
  \centering
  \begin{tabular}{p{0.6cm}|p{0.6cm}|p{0.6cm}|p{0.6cm}}
  \hline\hline
  $m$    & 2 & 3 & 4  \\ \hline
  $\eta$ & 30 & 30 & 30  \\ \hline
  $\# S$ & 12 & 20 & 25 \\ 
  \hline\hline
  \end{tabular}
\end{minipage}
\hspace{2cm}
\begin{minipage}[t]{0.3\textwidth}
  \begin{tabular}{p{0.6cm}|p{0.6cm}|p{0.6cm}}
    \hline\hline
   $m$ & 2 & 3  \\ \hline
   $\eta$ & 40 & 40   \\ \hline
   $\# S$ & 20 & 40  \\ 
    \hline\hline
  \end{tabular}
\end{minipage}
\caption{The $\# S$ used in 2D and 3D examples.} 
\label{tab_patch}
\end{table}

\begin{figure}[htbp]
  \centering
  \includegraphics[width=0.30\textwidth]{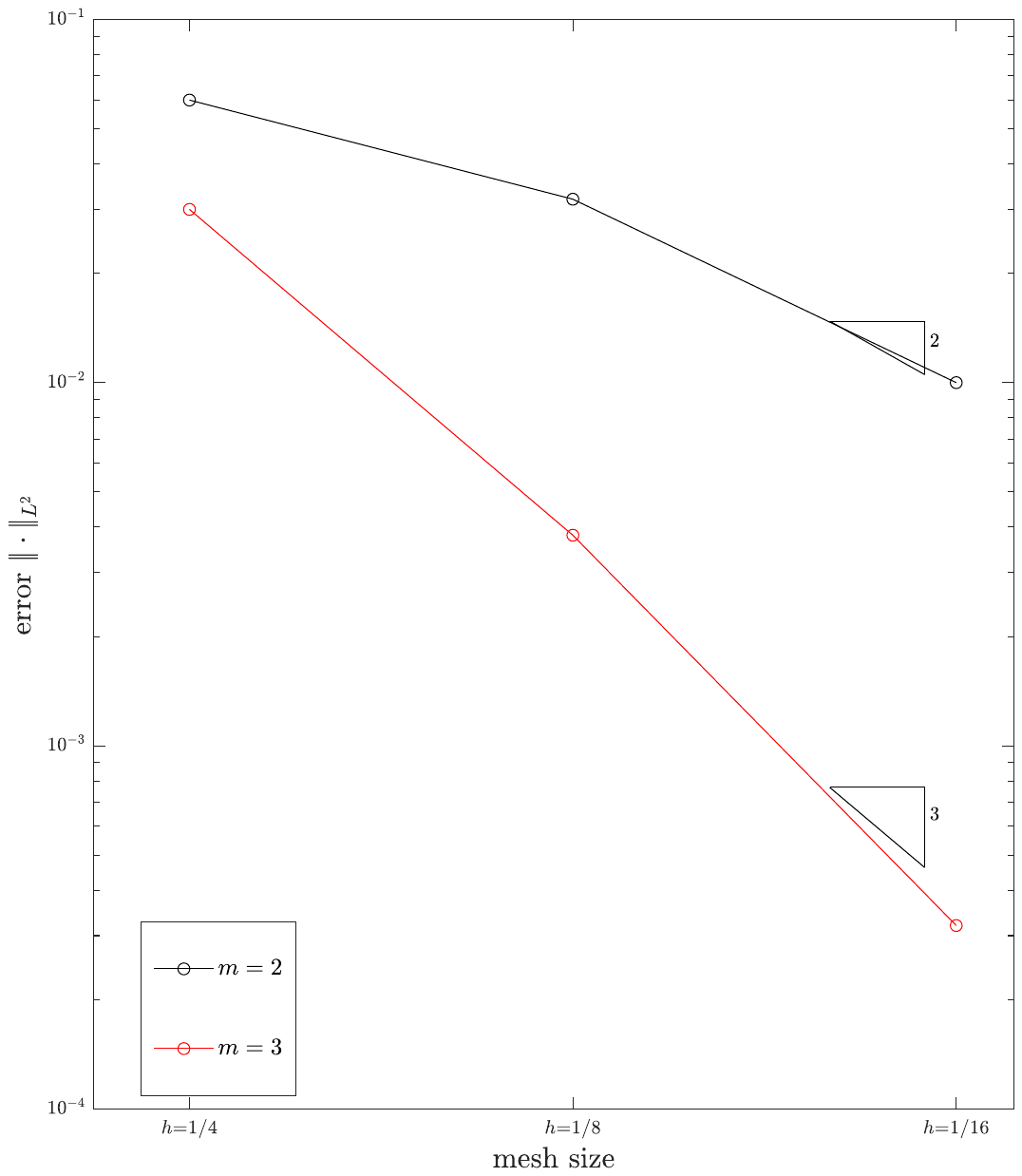}
  \hspace{30pt}
  \includegraphics[width=0.30\textwidth]{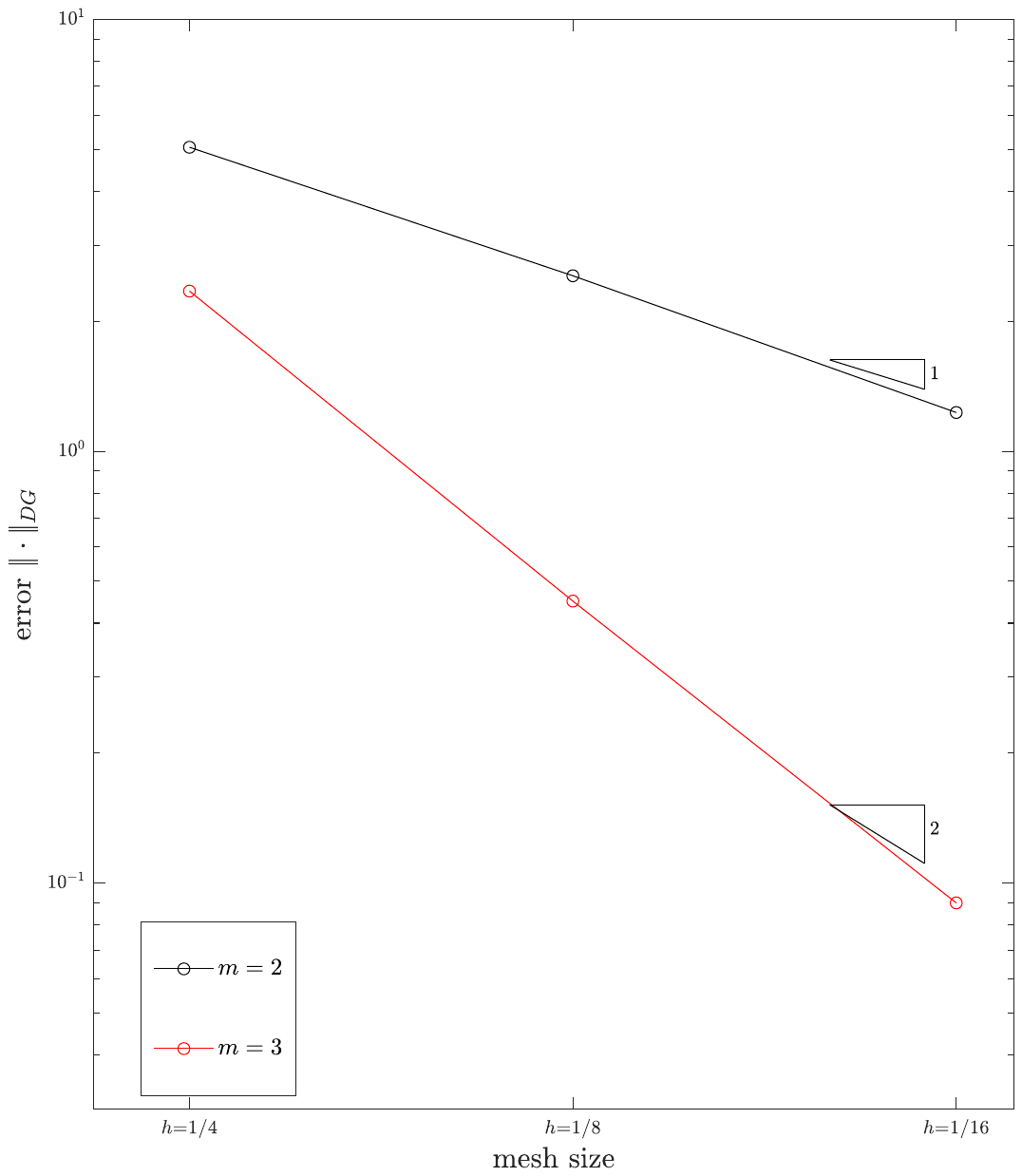}
  \caption{The convergence histories under the $\| \cdot \|_{L^2(\Omega)}$
  (left) and the  $\DGnorm{\cdot}$(right) in Example 2.}
  \label{fig_ex5err}
\end{figure}

\noindent \textbf{Example 2.}
Here we solve a 3D problem on $\Omega = (0,1)^3$.
We select the exact solution as
\begin{displaymath}
  \bm{u}(x, y, z) = \begin{bmatrix}
    \sin(\pi y)\sin(\pi z)\\
    \sin(\pi z)\sin(\pi x) \\
    \sin(\pi x)\sin(\pi y)\\    
  \end{bmatrix},
\end{displaymath}
We discretize the problem on successively refined meshes
with $h=1/4, 1/8, 1/16$. The convergence order is shown in
Fig.~\ref{fig_ex5err}, which confirms our theoretical result.

\section{Conclusion}
\label{sec_conclusion}
In this paper, we introduce an arbitrary order discontinuous Galerkin finite element method
to address the quad-curl problem. The discretization is based on a mixed method approach.
The approximation space is built using a patch reconstruction operator,
ensuring that the number of degrees of freedom remains unaffected by 
the approximation order. We establish optimal convergence in the energy norm 
and suboptimal convergence in the $L^2$ norm. 
Furthermore, we conduct numerical experiments in both two and three dimensions to validate our theoretical findings.

\bibliographystyle{amsplain}
\bibliography{../ref}

\providecommand{\bysame}{\leavevmode\hbox to3em{\hrulefill}\thinspace}
\providecommand{\MR}{\relax\ifhmode\unskip\space\fi MR }
\providecommand{\MRhref}[2]{%
  \href{http://www.ams.org/mathscinet-getitem?mr=#1}{#2}
}
\providecommand{\href}[2]{#2}
\begin{thebibliography}{10}

\bibitem{Ben2001finite}
N.~Ben~Salah, A.~Soulaimani, and W.~G. Habashi, \emph{A finite element method
  for magnetohydrodynamics}, Comput. Methods Appl. Mech. Engrg. \textbf{190}
  (2001), no.~43-44, 5867--5892. \MR{1848902}

\bibitem{Brenner2007mathematical}
S.~C. Brenner and L.~R. Scott, \emph{{The Mathematical Theory of Finite Element
  Methods}}, third ed., Texts in Applied Mathematics, vol.~15, Springer, New
  York, 2008.

\bibitem{Cakoni2010inverse}
F.~Cakoni, D.~Colton, P.~Monk, and J.~Sun, \emph{The inverse electromagnetic
  scattering problem for anisotropic media}, Inverse Problems \textbf{26}
  (2010), no.~7, 074004, 14. \MR{2644031}

\bibitem{Cakoni2007variational}
F.~Cakoni and H.~Haddar, \emph{A variational approach for the solution of the
  electromagnetic interior transmission problem for anisotropic media}, Inverse
  Probl. Imaging \textbf{1} (2007), no.~3, 443--456. \MR{2308973}

\bibitem{geuzaine2009gmsh}
C.~Geuzaine and J.~F. Remacle, \emph{Gmsh: {A} 3-{D} finite element mesh
  generator with built-in pre- and post-processing facilities}, Internat. J.
  Numer. Methods Engrg. \textbf{79} (2009), no.~11, 1309--1331.

\bibitem{Guermond2007MHD}
J.-L. Guermond, R.~Laguerre, J.~L\'{e}orat, and C.~Nore, \emph{An interior
  penalty {G}alerkin method for the {MHD} equations in heterogeneous domains},
  J. Comput. Phys. \textbf{221} (2007), no.~1, 349--369. \MR{2290574}

\bibitem{Han2023hp}
J.~Han and Z.~Zhang, \emph{An {$hp$}-version interior penalty discontinuous
  {G}alerkin method for the quad-curl eigenvalue problem}, BIT \textbf{63}
  (2023), no.~4, Paper No. 56, 29. \MR{4666352}

\bibitem{zhimin2023}
Jiayu Han and Zhimin Zhang, \emph{An hp-version interior penalty discontinuous
  galerkin method for the quad-curl eigenvalue problem}, BIT Numerical
  Mathematics \textbf{63} (2023), article number 56.

\bibitem{Hong2012discontinuous}
Q.~Hong, J.~Hu, S.~Shu, and J.~Xu, \emph{A discontinuous {G}alerkin method for
  the fourth-order curl problem}, J. Comput. Math. \textbf{30} (2012), no.~6,
  565--578. \MR{3041683}

\bibitem{Hu2020simple}
K.~Hu, Q.~Zhang, and Z.~Zhang, \emph{Simple curl-curl-conforming finite
  elements in two dimensions}, SIAM J. Sci. Comput. \textbf{42} (2020), no.~6,
  A3859--A3877. \MR{4186537}

\bibitem{Li2023curl}
R.~Li, Q.~Liu, and F.~Yang, \emph{A reconstructed discontinuous approximation
  on unfitted meshes to {$H({\rm curl})$} and {$H({\rm div})$} interface
  problems}, Comput. Methods Appl. Mech. Engrg. \textbf{403} (2023), no.~part
  A, Paper No. 115723, 27.

\bibitem{Li2016discontinuous}
R.~Li, P.~Ming, Z.~Sun, and Z.~Yang, \emph{An arbitrary-order discontinuous
  {G}alerkin method with one unknown per element}, J. Sci. Comput. \textbf{80}
  (2019), no.~1, 268--288.

\bibitem{Li2012efficient}
R.~Li, P.~Ming, and F.~Tang, \emph{An efficient high order heterogeneous
  multiscale method for elliptic problems}, Multiscale Model. Simul.
  \textbf{10} (2012), no.~1, 259--283.

\bibitem{Li2019reconstructed}
R.~Li and F.~Yang, \emph{A reconstructed discontinuous approximation to
  {M}onge-{A}mp\`ere equation in least square formulation}, Adv. Appl. Math.
  Mech. \textbf{15} (2023), no.~5, 1109--1141. \MR{4613677}

\bibitem{Monk2003finite}
P.~Monk, \emph{Finite element methods for {M}axwell's equations}, Numerical
  Mathematics and Scientific Computation, Oxford University Press, New York,
  2003.

\bibitem{Monk2012finite}
P.~Monk and J.~Sun, \emph{Finite element methods for {M}axwell's transmission
  eigenvalues}, SIAM J. Sci. Comput. \textbf{34} (2012), no.~3, B247--B264.
  \MR{2970278}

\bibitem{Powell1981approximation}
M.~J.~D. Powell, \emph{Approximation theory and methods}, Cambridge University
  Press, Cambridge-New York, 1981.

\bibitem{Sun2016mixed}
J.~Sun, \emph{A mixed {FEM} for the quad-curl eigenvalue problem}, Numer. Math.
  \textbf{132} (2016), no.~1, 185--200. \MR{3439219}

\bibitem{Wang2019new}
C.~Wang, Z.~Sun, and J.~Cui, \emph{A new error analysis of a mixed finite
  element method for the quad-curl problem}, Appl. Math. Comput. \textbf{349}
  (2019), 23--38. \MR{3894188}

\bibitem{wgcurl}
Chunmei Wang, Wang Junping, and Shangyou Zhang, \emph{Weak galerkin finite
  element methods for quad-curl problems}, Journal of Computational and Applied
  Mathematics \textbf{428} (2023), 115186.

\bibitem{Zhang2019curlcurl}
Q.~Zhang, L.~Wang, and Z.~Zhang, \emph{{$H({\rm curl}^2)$}-conforming finite
  elements in 2 dimensions and applications to the quad-curl problem}, SIAM J.
  Sci. Comput. \textbf{41} (2019), no.~3, A1527--A1547. \MR{3949709}

\bibitem{zhangqian2020}
Qian Zhang, \emph{A family of curl-curl conforming finite elements on
  tetrahedral meshes}, CSIAM Transactions on Applied Mathematics \textbf{1}
  (2020), 639--663.

\bibitem{Zhang2018mixed}
S.~Zhang, \emph{Mixed schemes for quad-curl equations}, ESAIM Math. Model.
  Numer. Anal. \textbf{52} (2018), no.~1, 147--161. \MR{3808156}

\bibitem{zheng2011nonconforming}
B.~Zheng, Q.~Hu, and J.~Xu, \emph{A nonconforming finite element method for
  fourth order curl equations in {$\Bbb{R}^{3}$}}, Math. Comp. \textbf{80}
  (2011), no.~276, 1871--1886. \MR{2813342}
  
\bibitem{2019cc}
Jiguang Sun, Qian Zhang, and Zhimin Zhang, \emph{A curl-conforming weak
  galerkin method for the quad-curl problem}, BIT Numerical Mathematics
  \textbf{59} (2019).
\end{thebibliography}

\end{document}